\documentclass[A4, reqno]{amsart}
\usepackage[dvips]{graphicx}
\usepackage{euscript}
\usepackage{amscd,xy}
\usepackage{amssymb,latexsym}

\usepackage{amssymb}
\xyoption{all}
\newtheorem{thm}{Theorem}[section]
\newtheorem{dfn}[thm]{Definition}
\newtheorem{prop}[thm]{Proposition}
\newtheorem{cor}[thm]{Corollary}
\newtheorem{lma}[thm]{Lemma}
\newtheorem{rem}[thm]{Remark}

\newenvironment{pf}{\begin{proof}}{\end{proof}}

\numberwithin{equation}{section}

\newcommand{\R}{{\mathbb{R}}}

\newcommand{\C}{{\mathbb{C}}}
\newcommand{\Q}{{\mathbb{Q}}}
\renewcommand{\P}{{\mathbb{P}}}

\begin{document}

\title[Topological quandles and invariants of links]
{Topological quandles and invariants of links}

\author{Ryszard L. Rubinsztein\\ \\ 2005-08-26}
\address{Dept. of Mathematics, Uppsala University, Box 480,
Se-751 06 Uppsala, Sweden}
\email{ryszard{\@@}math.uu.se}

\begin{abstract}
We introduce a notion of topological quandle. Given a topological quandle $Q$ we associate to every classical link $L$ in $\R ^3$ an invariant $J_Q(L)$ 
 which is a topological space (defined up to a homeomorphism). The space  $J_Q(L)$  can be interpreted as a space of colourings of a diagram of the link   $L$ 
with colours from the quandle $Q$. 
\end{abstract}

\subjclass[2000]{Primary: 57M27, 57M25 ;
Secondary: 20F36 }

\keywords{Topological quandle, invariant spaces of links, spaces of colourings of links. }

\maketitle

\section{Introduction}

\bigskip
  
Quandles (and related objects: automorphic sets and racks) have been investigated in connection with knots at least from the beginning of 1980-ies, \cite{J1}, \cite{M1}, 
\cite{B1}, \cite{F1}. They are (discrete) sets equipped with a binary operation satisfying certain identities.

In this paper  we introduce a notion of topological quandles. These are quandles equipped with topology so that the binary operation is continuous, see Definition 2.1.

Given a topological quandle $\, Q\,$ we define an invariant of (the ambient isotopy classes of) oriented classical links in $\, \R ^3\,$. 
For such a link $\, L\,$ the invariant $\, J_Q(L)\,$ is a topological space 
(well-defined up to a homeomorphism). The space  $\, J_Q(L)\,$ is defined as the space of fix-points of an action of a braid the closure of which represents the 
link $\, L\,$.  There are two other possible interpretations of  $\, J_Q(L)\,$ : as  a space of colourings of a diagram of   $\, L\,$ with colours from $\, Q\,$ 
and as the space of quandle homomorphisms $\, Hom_q(X(L),\, Q)\,$ from the universal (discrete) quandle $\, X(L)\,$ of the link  $\, L\,$ (as defined by Joyce and Matveev) 
to the topological quandle $\, Q\,$.

\medskip

The paper is organized as follows. In Section 2 we define topological quandles and give examples of such quandles. In Section 3 we recall the definition of an action of 
braid groups on products of quandles. Section 4 contains the definition of the invariant space   $\, J_Q(L)\,$ of a link and two other interpretations of this space. 
The main result of  the section (and of the paper) is Theorem 4.1 stating that for every topological quandle $\,Q\,$ the space  $\, J_Q(L)\,$ depends only 
on the isotopy class of the oriented link  $\, L\,$. 
In Section 5 we calculate the invariant space  $\, J_Q(L)\,$ for the Hopf link, for the trefoil knot and for the figure-eight knot 
in the case when the quandle $\, Q\,$ is the 2-dimensional sphere $\, S^2\,$  with the quandle structure given by $\, S^2\,$  being a symmetric Riemannian manifold. 
Finally, in Section 6 we calculate the invariant space  $\, J_Q(L)\,$ for the trefoil knot in the case when the quandle $\, Q\,$ is a conjugacy class in the Lie group 
$\, SL(2,\, \C\, )\,$.
\medskip

I am greatly indebted to Oleg Viro for many enlightening discussions on the subject of this paper. I am also grateful to Volodia Mazorchuk for useful comments 
and for his help 
with the graphics of the paper. 
\section{Topological quandles}

\medskip 
\begin{dfn}\label{dfn 2.1}
Let $X$ be a topological space equipped with a continuous map $\, \mu :X\times X\longrightarrow X\,$, denoted by $\, \mu (a, b) = a\star b ,\, $ such that for every 
$\, b\in X\,$ the mapping $\, a\mapsto a\star b\,$ is a homeomorphism of $\,X\,$.
The space $\, X\,$ (together with the map $\, \mu \,$) is called a {\bf topological quandle} if it satisfies

(i) \,\, $\, (a \star b) \star c = (a\star c)\star (b\star c)\,$,

(ii) \,\, $\, a\star a = a\,$,

\noindent
for all $\, a, b, c \in X\,$.
\end{dfn}

 A continuous map $\, \varphi : X\longrightarrow Y\,$ between topological quandles $\, X\,$ and $\, Y\,$ is called a {\bf quandle homomorphism} 
if $\, \varphi (a \star b) = \varphi (a) \star \varphi (b)\,$ holds for all $\, a,b \in X\,$.

A subspace $\, Z\,$ of a topological quandle $\, X\,$ is called a {\bf subquandle} of $\, X\,$ if $\, a\star b\in Z\,$  for all $\, a,b\in Z\,$ and, 
for every $\, b\in Z\,$, the mapping $\, z\mapsto z\star b\,$ is a homeomorphism of $\, Z\,$ onto itself. 
\bigskip

\noindent 
{\bf Example 2.1: Symmetric manifolds.} Let $M$ be a symmetric manifold  i.e. $ M\,$ is a Riemannian manifold such that each point $\, b\in M\,$ 
is an isolated fixed point of an involutive isometry $\, i_b:M\longrightarrow M\,$, (see \cite{H1}). 
Then the operation $\, a\star b = i_b(a)\,$ provides $M$ with a structure of a topological quandle. Riemannian isometries of $\, M\,$ are 
quandle endomorphisms of this quandle. 

\medskip

\noindent
{\bf Example 2.2: Conjugation quandles.} Let $G$ be a topological group. For $g, h \in G$ define $\, h\star g: = g^{-1}hg\,$. This turns $G$ into a topological quandle 
called 
the conjugation quandle of $G$. Continuous group homomorphisms of topological groups are quandle homomorphisms of the corresponding conjugation quandles.

\medskip

\noindent
{\bf Example 2.3: Alexander quandles.} Let $\, \sigma : G\longrightarrow G\, $ be a (continuous) automorphism of a topological group $G$. The operation 
$\, h\star g: = \sigma (g)^{-1} \sigma (h) g\, $ turns $G$ into 
a topological quandle. (Proof of that is a matter of a simple check and left to the reader). Example 2.2 is a special case of it with $\,\sigma = id\,$.

\medskip

\noindent
{\bf Example 2.4: Anti-Alexander quandles.} Let $\, \tau : G\longrightarrow G\,$ be a (continuous) anti-automorphism of a topological group $G$. We can compose it with  
the mapping $\, g\mapsto g^{-1}\,$, 
obtaining an automorphism of $G$, and then apply the construction of Example 2.3 to obtain a quandle structure on $G$. Explicitly that gives 
 $\, h\star g: = \tau (g) \tau (h)^{-1}g\,$. 
Again, Example 2.2  is a special case of this with $\, \tau (g) = g^{-1}\,$.

\medskip

\noindent
{\bf Example 2.5: Conjugacy subquandles.} Let again $G$ be a topological group and $\, X\subset G\,$ be a subset which is a union of some cojugacy classes of $G$. 
The restriction of 
the operation from Example 2.2  to $\, X\,$ makes it a topological quandle. It is a subquandle of the conjugacy quandle of $\, G\,$.

\medskip

\noindent
{\bf Example 2.6:  Spheres.} Let $\, q\,$ be a complex number of modulus $1$. Let $\, V\,$ be a comlex vector space of finite dimension equipped with 
a hermitian inner product $\, < \cdot , \cdot >\,$. For $\, a \in V, \, a\ne 0,\,$ define a unitary isometry $\, i^q_a:V\longrightarrow V\,$ by 
\begin{displaymath} 
i^q_a (v) = q\, v + (1-q)\,\, \frac {<v , a>}{<a,a>} \,\,\, a \quad .
\end{displaymath}    
($\, i^q_a\,$ is the identity on the complex line spanned by $a$ and the multiplication by $q$ on the orthogonal complement of the line). Let $S$ be the unit sphere in $V$. 
For $\, a,b\in S\,$ define $\, a\star b : = i^q_b (a)\, $. That gives $S$ a structure of a topological quandle.

\begin{proof}
Let $\, V = l_a \oplus l_a^{\perp}\,$ be the splitting of $V$ into the orthogonal direct sum of the line $l_a$ spanned by $a$ and of its orthogonal complement. 
For a vector $v\in V$,
the corresponding splitting has the form $\, v=v_1 \oplus v_2\,$ with
\begin{displaymath}
v_1=v- \frac {<v , a>}{<a,a>} \,\,\, a  \qquad \text{and} \qquad  v_2= \frac {<v , a>}{<a,a>} \,\,\, a \quad .
\end{displaymath}  
Therefore $\, i^q_a (v) =  i^q_a ( v_1 \oplus v_2 ) = q \, v_1 \oplus v_2\,$ which shows that $\,  i^q_a\,$ is a unitary isometry (as $|q|=1$) and that  
$\, a\star a =  i^q_a (a) = a\,$ proving the property $(ii)$ of Definition 2.1.

The property $(i)$ follows from the fact that
\begin{equation}\label{(2.1)}
 i^q_b (\, i^q_a (v)\, ) =  i^q_{ i^q_b (a)}(\, i^q_b (v)\, )  % \qquad \qquad \qquad  (2.1)
\end{equation}
for all $\, a, b, v\in S\, $. This identity follows directly from the definition of $\,  i^q_b\,$ (and the fact that  $\,  i^q_b\,$ is a unitary isometry) and 
is left to the reader.. 
%To prove $(i)$ we have to show that if $\, a, b, v\in S\, $ then  $\,  i^q_b (\, i^q_a (v)\, ) =  i^q_{ i^q_b (a)}(\, i^q_b (v)\, )\,$.  Since $\, <a,a>=<b,b>=1\,$ we have 
%\begin{equation*}
%\begin{split}
% i^q_b (\, i^q_a (v)\, )& =  i^q_b (\,  q \, v  + (1-q) <v, a> a\, ) =\\ & = q^2v+q(1-q) <v, a> a +\\
%                        & \qquad \qquad +  (1-q) < q  v  + (1-q) <v, a> a \,\, ,\, b > \, b  = \\
%                        & =  q^2v+q(1-q)( <v, a> a \, + <v,b>\, b ) +\\ & \qquad \qquad + (1-q)^2 <v,a><a,b>\, b \quad ,
%\end{split}
%\end{equation*}
%while
%\begin{equation*}
%\begin{split}
% i^q_{ i^q_b (a)} (\, i^q_b (v)\, )& =  q\, i^q_b (v)\, + (1-q)\, < \, i^q_b (v)\, ,\, i^q_b (a)\, >\,  i^q_b (a)\, =\\
%                                   & =  q\, i^q_b (v)\, + (1-q)\, < \, v\, ,\, a\, >\,  i^q_b (a)\, = \\
%                                   & =  q\, (\, qv + (1-q)<v,b>\, b\, ) + \\ 
%                                   & \qquad\qquad + (1-q) <v,a> (\, qa + (1-q)<a,b>\, b\, ) = \\
%                                   & = q^2v + q(1-q)(<v,b>\, b\, + <v,a>\, a\, ) +\\ & \qquad \qquad  +  (1-q)^2 <v,a><a,b>\, b\quad ,
%\end{split}
%\end{equation*}
%which proves $(i]$.
\end{proof}
 
\medskip

\noindent
{\bf Example 2.7: Projective spaces.} Let $\, V\,$ and $\, q\,$  be as in Example 2.6 and let $\, \mathbb{P} (V) \,$ be the projective space of complex lines in  $\, V\,$. 
The quandle structure on 
the sphere $\, S\,$ in $\, V\,$ described in Example 2.6
induces a topological quandle stucture on  $\, {\P} (V) \,$.  If $\, l_a, \, l_b\,$ are lines in $\,V\,$ spanned by vectors $a$ resp. $b$ then 
$\, l_a \star l_b =  i^q_b (\, l_a\,)\,$. 

There are corresponding quandle structures for real vector spaces    $V$  and  $q=-1$, but these are special cases of Example 2.1.  

\medskip

\noindent
{\bf Example 2.8: Grassmannians.} Let $\, V\,$ and $\, q\,$  be again as in Example 2.6 and  $\, k\,$ be an integer 
such that $\, 1\le k \le \dim V\,$. Denote by $\, Gr_k(V)\,$  
the grassmannian of $k$-dimensional linear subspaces of $\, V\,$. For any subspace $\, U\in  Gr_k(V)\,$ let $\, j^q_U : V\longrightarrow V\,$ be the unitary isometry 
defined by 
\begin{displaymath} 
j^q_U (v) = q\, v + (1-q)\,\,\sum\limits _{i=1}^k\, <v , a_i> \,\, a_i \quad ,
\end{displaymath} 
where $\, a_1, ... , a_k\,$ is an orthonormal  basis of $\, U\,$. In other words, $\, j^q_U\,$ is the identity on $\, U\,$ and the  multiplication by $\, q\,$ on 
the orthogonal complement of $\, U\,$. For two subspaces $\, U, W \in Gr_k(V)\,$   one checks directly the identity
\begin{equation}
         j^q_W   ( j^q_U (v)) =  j^q_{ j^q_W (U)}   ( j^q_W (v))
\end{equation}
for all $\, v\in V\,$.

Now, for  $\, U, W \in Gr_k(V)\,$ define $\, U\star W: =  j^q_W (U)\,$. This operation gives   $\, Gr_k(V)\,$ a structure of a topological quandle. That is a consequence 
of the identity (2.2). 

Again, there are  corresponding quandle structures for real vector spaces    $V$  and  $q=-1$.

\bigskip
\noindent
{\bf Remark}: The notion of a quandle  as well as related notions of {\it kei}, of {\it authomorphic sets} and of {\it racks} appeared  in the mathematical literature 
on many occasions, see \cite{B1}, \cite{F1}, \cite{J1}, \cite{K1}, \cite{K2}, \cite{M1}, \cite{T1}.  However, in all these cases automorphic sets, quandles and racks were 
treated 
in the category of sets. Symmetric manifolds (Example 2.1) were the source of the definition of {\it kei} by M.Takasaki, \cite{T1}, in 1942. They were 
also mentioned in \cite{B1} and 
\cite{K2}, yet no use was made of their topological structure.

A good reference for a review of definitions, properties and applications of keis, quandles and racks (as sets) is \cite{K1}. 

In 1982,  independently, D.Joyce, \cite{J1}, and S.Matveev, \cite{M1}, introduced  the notion of a quandle of a knot. (Matveev used the name {\it distributive groupoid} 
instead of {\it  quandle}.)  Their quandle is again a discrete object 
(a set with an operation) and 
a complete invariant of oriented knots (up to ambient isotopy).  

\bigskip

We define also actions of topological quandles on spaces.

\medskip

\begin{dfn}\label{dfn 2.2}
Let $\, X\,$ be a topological quandle , $\, M\,$ be  a topological space and  
$\, \varphi : M \times X \longrightarrow M\,$ be a continuous map. We denote $\, \varphi ( m, a ) = m \star a\,$.

The map $\, \varphi : M \times X \longrightarrow M\,$ is an action of the quandle  $\, X\,$ on the space $\, M\,$ if
\begin{displaymath}
( m \star a ) \star b = (m \star b ) \star ( a \star b )
\end{displaymath}
holds for all $\, a, b \in X\,$ and $\, m\in M\,$.
\end{dfn}

\medskip

\noindent
{\bf Example 2.9}:   If $\, G\,$ is a topological group acting from the right on a topological space $\, M\,$ with 
$\,\, \phi : M\times G \longrightarrow M\,$ being the the action then the same map $\, \phi \,$ is also an action of the 
conjugation quandle of $\, G\,$ on the space $\, M\,$.
Indeed, denoting the ordinary action of the group element $a \in G$ on a point $m\in M$ by $ma$,  
we have $\, ( m \star a ) \star b = (ma)b = m(ab) = (mb)(b^{-1}ab) = (m \star b ) \star ( a \star b)\,$.

\medskip

\noindent
{\bf Example 2.10}: Let $\, V\,$ be a finite dimensional hermitian vector space and $\ q\,$ a complex number, $\, |q|=1\,$, as in Examples 2.6 and 2.7. 
Let $\, \P (V)\,$ be the projective space of lines in $\, V\,$ equipped with the topological quandle structure from Example 2.7. Let $\, Fl (V)\,$ be a space of flags 
(complete or partial of some kind) in $\, V\,$. Then there is an action of the quandle   $\, \P (V)\,$ on the space   $\, Fl (V)\,$ given by 
$\, f \star l_a : = i^q_a (f)\,$ for any $\, f\in Fl(V)\,$ and $\, l_a\in \P (V)\,$. That this is indeed an action of a quandle on a space follows directly from 
the formula (2.1).   

\medskip

\section{Braids and quandles}

\medskip

Let $\, B_n\,$ be the classical braid group on $n$ strands and let $\, Q\,$ be a quandle. It has been observed in \cite{B1} that the braid group   $\, B_n\,$ 
acts (from the right) on the product $\, Q  \times Q \times  ... \times Q\,$   of n copies of  $\, Q\,$. We shall now describe this action.

Let us recall that the group $\, B_n\,$ is generated by $\, n-1\,$ elementary braids $\, \sigma _i, \,\, i=1, ..., n-1,\,$ subject to relations 
\begin{equation}
\begin{split}
&\sigma _i \, \sigma _j \, = \sigma _j \, \sigma _i \, \qquad \text{if} \quad |i-j|>1\quad ,\\
&\sigma _i \, \sigma _{i+1} \,  \sigma _i =  \sigma _{i+1} \, \sigma _i \,  \sigma _{i+1} \, \qquad \text{for} \quad i = 1, ... , n-2.
\end{split}
\end{equation}
The right action of the elementary braid $\, \sigma _i\,$ on  the product $\, Q  \times Q \times  ... \times Q\,$   of n copies of  $\, Q\,$ is defined by 
%\begin{displaymath}
%\sigma _i (x_1,...,x_{i-1}, x_i, x_{i+1},x_{i+2},...,x_n) =  (x_1,...,x_{i-1}, x_{i+1}, x_i\star x_{i+1}, x_{i+2},...,x_n).
%\end{displaymath}
\begin{equation}
(\sigma _i (x))_j = \left\{
\begin{array}{lcl}
x_j& \text{if} & j\ne i, \, i+1,\\
x_{i+1} & \text{if} & j = i,\\
x_i\star x_{i+1} & \text{if} & j = i+1,
\end{array} 
\right.
\end{equation}
where $\, x_j \in Q\,$ is the $j$-th coordinate of a point $\, x\in   Q  \times Q \times  ... \times Q\,$.

One checks that the right action of the elementary braids $\, \sigma _i\,$  so defined satisfies the relations (3.1). The first one is obvious and the second one 
is a matter of 
a direct check. Hence, the formulas (3.2) define a right action of the braid group $\, B_n\,$ on  the product  $\, Q  \times Q \times  ... \times Q\,$ of $n$ copies 
of the quandle  $\, Q\,$. If $\,Q\,$ is a topological quandle the action is continuous.  

\smallskip
\noindent
{\bf Remark 3.1 :} If $\, \varphi : Q_1 \longrightarrow Q_2\,$ is a quandle homomorhism then 
the product map $\, \prod\limits ^n \,  \varphi :\prod\limits ^n \, Q_1 \longrightarrow  \prod\limits ^n \,   Q_2\,$ commutes with the action of the braid group 
$\, B_n\,$ on $\, \prod\limits ^n \, Q_1\,$ and  $\, \prod\limits ^n \, Q_2\,$.
\medskip

\section{Invariants of links from topological quandles}

\medskip

In this section, given a topological quandle Q, we construct an invariant of links. The invariant will associate to every oriented link $\, L\,$ in $\, \R ^3\,$  
a topological space  $\, J_Q(L)\,$ (well-defined up to a homeomorphism). 

Let $\, Q\,$ be a topological quandle fixed for the rest of this section.

Any oriented link $\, L\,$ in $\, \R ^3\,$ can be described as a closure of a braid $\, \sigma\,$ on $n$ strands, $\,\sigma\in B_n\,$, for some $n\ge 1$. 
According to Section 2 
there is the action of $\, B_n\,$ on the  product  $\, Q  \times Q \times  ... \times Q\,$ of $n$ copies 
of the quandle  $\, Q\,$. This action gives us, in particular, a continuous map 
\begin{displaymath}
\sigma :  Q  \times Q \times  ... \times Q \longrightarrow  Q  \times Q \times  ... \times Q \quad .
\end{displaymath}
We define $\, J_Q(L)\,$ to be the space of fix-points of  the map $\, \sigma\,$. Thus  $\, J_Q(L)\,$ is a subspace of $\, Q  \times Q \times  ... \times Q\,$ .

\medskip
\begin{thm}
For every topological quandle $\, Q\,$  the topological space  $\, J_Q(L)\,$ is an invariant of the oriented link $\, L\,$ i.e. it depends (up to a homeomorphism) only on 
the isotopy class of $\, L\,$ in $\, \R ^3\,$.   
\end{thm}

\begin{pf}
To see that the space $\, J_Q(L)\,$ depends only on the link   $\, L\,$ and not on a particular choice of the braid  $\, \sigma\,$ representing it, we have to show, 
according to the theorem of Markov, \cite{B2},  that  

(i) any braid conjugate to   $\, \sigma\,$  in $\, B_n\,$ has the space of fix-points homeomorphic 

\quad\,  to that of   
$\, \sigma\,$  itself, and  

(ii) the space of fix-points of $\, \sigma\,$ is homeomorphic to that of the braid   $\, \sigma \,  \sigma _n^{\pm 1} \in $

$\quad\,\,\, B_{n+1}   \,$. 

\noindent
(In this last statement we use the standard inclusion of $\, B_n\,$ into   $\, B_{n+1}\,$.)

The property $(i)$ is obvoius since conjugate homeomorphisms have homeomorphic fix-points sets. We shall now prove the property $(ii)$. 

Let $\, \sigma \in B_n\,$ and let 
\begin{displaymath}
\sigma : \, \prod \limits ^n Q \longrightarrow  \prod \limits ^n Q   \qquad  (n \,\,\,  \text{factors of} \,\, \,  Q)
\end{displaymath}
and
\begin{displaymath}
\sigma _n \circ \sigma  :\,\, \prod \limits ^{n+1} Q \longrightarrow  \prod \limits ^{n+1} Q  \qquad  (n+1 \,\,\,\,  \text{factors of} \,\, \,  Q )
\end{displaymath}

\medskip

\noindent
be the corresponding maps (remember that our action is a right one). If $\, x=(x_1, ... , x_{n+1}) \in  \prod \limits ^{n+1} Q  \,$ is a fixed point 
of $\,\sigma _n \circ \sigma  \,$ we shall show that $\, x_n = x_{n+1}\,$ and that $\,\tilde x = (x_1, ... , x_n)\in  \prod \limits ^n Q \,$ is a fixed point of 
$\, \sigma\,$. 

Indeed, let $\, \sigma (x_1, ... ,x_n) = ( y_1, ... ,y_n) \in  \prod \limits ^n Q \,$. Then 
\begin{equation*}
\begin{split}
(x_1, ... ,x_n,  x_{n+1}) &= ( \sigma _n \circ \sigma ) (x_1, ... , x_n, x_{n+1})=  \sigma _n ( \, \sigma (x_1, ... ,x_n),  x_{n+1})=\\
& =  \sigma _n ( \, y_1, ... , y_n, x_{n+1}) = \\
& = (y_1, ... , y_{n-1},  x_{n+1}, y_n \star  x_{n+1}).
\end{split}
\end{equation*} 

It follows that 
\begin{equation}
\begin{split}
& y_i = x_i \qquad \text{for} \quad i=1, ... , n-1,\\
& x_{n+1} = x_n,\\
& y_n \star  x_{n+1} =  x_{n+1}.
\end{split}
\end{equation} 

Since, according to the definition of a quandle, the mapping $\, a \mapsto  a\star  x_{n+1} \,$ is a homeomorphism of $\, Q\,$ onto itself  and since 
$\,  x_{n+1}  \star  x_{n+1} =  x_{n+1}\,$ , the last equation of (4.1) implies that $\, y_n = x_{n+1}\,$. Together with the second equation of (4.1) this gives  
\begin{displaymath}
 \sigma (x_1, ... ,x_n) = ( y_1, ... ,y_{n-1}, y_n) = (x_1, ..., x_{n-1}, x_n )
\end{displaymath}
i.e. $\, \tilde x = (x_1, ... ,x_n) \in  \prod \limits ^n Q \,$ is a fixed point of 
$\, \sigma\,$. By the second equation of (4.1) we have also $\,  x_{n+1} = x_n\,$. 

On the other hand, it is immediate that if $\, \tilde x = ( x_1, ... ,x_n)\in \prod \limits ^n Q  \,$ is a fixed point of $\, \sigma \,$ then 
$\, x = (x_1, ... ,x_n,  x_n) \in    \prod \limits ^{n+1} Q \, $ is a fixed point of $\,  \sigma _n \circ \sigma\,$ .  

All this shows that the projection 
onto the first n factors $\,   \prod \limits ^{n+1} Q \longrightarrow \prod \limits ^n Q \,$  induces a homeomorphism from the space of fix-points 
of  $\,  \sigma _n \circ \sigma\,$ onto the space of fix-points of  $\, \sigma \,$. 

One proves in the same way that the fix-point sets of  $\, \sigma \,$ and of $\, \sigma _n^{-1} \circ \sigma \,$ are homeomorphic. 
That finishes the proof of the property $(ii)$ and, hence, shows that the space $\, J_Q(L)\,$ is well-defined up to a homeomorphism .
\end{pf}

\medskip

{\bf Remark 4.1 :} Any quandle homomorphism $\, \varphi  : Q_1 \longrightarrow Q_2\,$ induces a continuous 
map $\, J_{\varphi }(L) : J_{Q_1}(L) \longrightarrow J_{Q_2}(L)\,$. Indeed, the map $\,  J_{\varphi }(L)\,$ is the restriction of $\, \prod\limits ^n \,  \varphi\,$ 
to $\, J_{Q_1}(L)\,$. The assignment $\, \varphi \longmapsto  J_{\varphi }(L)\,$ is functorial.

\bigskip

Below we shall give two other descriptions of $\,  J_Q(L)\,$.

\medskip

For the first one, let us describe the link $\, L\,$ by an oriented link diagram (in the oriented plane). For the moment we choose a link diagram given by a description 
of the link as 
a closure of a braid $\sigma \in B_n\,$ as above. Later on, it will be clear that this assumption on the kind of the link diagram is unnecessary and that any link 
diagram can be used for the purpose of descibing the space $\,  J_Q(L)\,$.

Now, let us consider the set of all colourings of the link diagram by elements of the quandle $\, Q\,$ defined as follows.
To every arc of the link diagram we associate an element of $\, Q\,$ subject only to the requirment that at any crossing of the diagram 
the following should be satisfied: if $\, a\in Q\,$ is associated to the directed {\it overcrossing} arc of the crossing,  $\, b\in Q\,$ is 
associated to the arc laying on the left-hand side of the overcrossing arc and $\, c\in Q\,$ is associated to the arc laying on the right-hand side 
of the overcrossing arc then 
\begin{equation} 
c =  b\star a  \quad . 
\end{equation}
See Figure 1.
\begin{displaymath}
\xymatrix{
a& c\\
b\ar@{-}[ru]|\hole  & a\ar[lu]
}\begin{array}{c}
\\ \\ \\ \qquad \qquad c=b\star a\\
\end{array}
\end{displaymath}

\centerline{Figure 1.}

\medskip

Here we adopt the convention that braids are oriented from the bottom to the top. The start- and endpoints of the strands of the braids are ordered from 
the left to the right.
 
Every such a  colouring of the diagram is uniquely determined by the sequence of $n$ elements  
$\, (x_1, ..., x_n) \in  \prod \limits ^n Q \,$, where $\, x_i\,$ is the element of $\,Q\,$  associated to the arc at the starting point of the $i$-th strand of the braid. 
Moving up along the braid the sequence  $\, (x_1, ..., x_n) \,$ determines the colours (elements) associated to all arcs of the non-closed braid. It follows that
the sequence of elements associated to the arcs at the endpoints (at the top) of the braid is equal to $\, \sigma (x_1, ..., x_n) \in  \prod \limits ^n Q \,$. To be able 
to colour the {\it closed} braid (and hence the diagram of the link) one needs that the start points and the end points of the braid with the same index are coloured by the 
same colour i.e. one needs the equality
\begin{displaymath}
\sigma (x_1, ..., x_n) = (x_1, ..., x_n) \quad .
\end{displaymath}
Hence, a sequence $\, (x_1, ..., x_n)  \in  \prod \limits ^n Q \,$ determines a colouring of the link diagram if and only if it is a fix-point of 
the action of $\sigma \,$ on $\,  \prod \limits ^n Q   \,$,  i.e. if and only if   $\, (x_1, ..., x_n)  \in J_Q(L) \,$. 
In that way the set of all colourings of the link diagram can be identified with $\, J_Q(L) \,$.  This description of the space   $\, J_Q(L) \,$ was suggested by Oleg Viro.

\smallskip

To get the third description of the invariant space  $\, J_Q(L) \,$, let us recall the definition of the (discrete) link quandle $\, X(L)\,$ 
of the link $\, L\,$ introduced by 
D.Joyce, \cite{J1},  and S.Matveev, \cite{M1}. We follow \cite{K1},Section 3, except that we use an opposite orientation of the plane to be in agreement with our previous 
definitions.  

Let us choose an oriented plane diagram $\, D\,$ of the link $\, L\,$. Let $\, E\,$ be the set of arcs of  $\, D\,$. Then the link quandle $\, X(L)\,$ is 
the discrete quandle 
generated by the set $\, E\,$ subject only to the relations coming from the crossings of the diagram in the following way: if $\, e\,$ is the {\it overcrossing} arc 
of the crossing, $\, f\,$ is the arc to the left of $\, e\,$ and $ \, g \,$ is the arc to the right of   $\, e\,$ then we get a relation 
\begin{equation}
g=f\star e \quad 
\end{equation}   
See Figure 2.\begin{displaymath}
\xymatrix@!=0.3pc{
&& \\
f\ar@{-}[rr]|\hole&  & g\\
&e\ar[uu]&
}\begin{array}{c}
\\ \\ \\ \\ \qquad \qquad g=f\star e\\
\end{array}
\end{displaymath}
\centerline{Figure 2.}

\medskip

\newpage
\noindent
The quandle  $\, X(L)\,$ so defined  depends only on the isotopy class of the link  $\, L\,$ and not on a particular choice of the diagram   $\, D\,$. 

Given a topological quandle $\, Q\,$, let us consider the set of all quandle homomorphisms $\, Hom_q( \, X(L),\, Q \, )\,$ equipped with the compact-open topology. 
If $\, D\,$ is 
an oriented plane diagram of the link $\, L\,$, it is clear from the definition of  $\, X(L)\,$ that every element $\, \varphi \in  Hom_q( \, X(L),\, Q \, )\,$ is given 
by associating an element of $\, Q\,$ to every arc of the diagram $\, D\,$ in such a way that the relations (4.2) are satisfied at every crossing of  $\, D\,$.
In other words, elements of $\, Hom_q( \, X(L),\, Q \, )\,$ can be identified with the colorings of the diagram  $\, D\,$. This observation gives a homeomorphism between the 
 space   $\, Hom_q( \, X(L),\, Q \, )\,$ and the space of colourings of the diagram $\,D\,$  Here the space of colourings is equipped with the relative topology of a subset 
of the product 
$\, \prod\limits ^m  Q\, , \,\, \, m$ being the number of arcs of $\, D\,$.    

All this holds for any oriented diagram  $\, D\,$ 
of the link  $\, L\,$. In particular, choosing a diagram coming from a description of  $\, L\,$ as a closure of a braid we get

\begin{cor}  
The spaces $\, J_Q(L) \,$ and  $\, Hom_q( \, X(L),\, Q \, )\,$ are homeomorphic.
\end{cor}
   
\medskip

\begin{cor}  
For any  plane diagram $\, D\,$ of the oriented link $\, L\,$  the space of colourings of $\, D\,$ with the elements of a topological quandle $\, Q\,$ is homeomorphic 
to the space 
$\, J_Q(L) \,$.
\end{cor}

\bigskip
\noindent
{\bf Remark 4.3.} If the quandle $\, Q\,$ satisfies an additional condition that
\begin{equation}
(b\star a)\star a = b 
\end{equation}
for all $\, a,b \in Q\,$ (i.e.  if $\, Q\,$ is a {\it kei}, \cite{K1}) then the invariant space $\, J_Q(L)\,$ does {\bf not} depend on the orientation of the link $\, L\,$. 
That is, perhaps, most easily seen when interpreting   $\, J_Q(L)\,$ as the space of colourings of a diagram $\, D\,$ of   $\, L\,$. In Figure 1, with the given 
orientation of $\, L\,$, we have the relation $\, c=b\star a\,$. If we change the orientation of  $\, L\,$, the orientation of the overcrossing changes and 
the corresponding relation is  $\, b=c\star a\,$. When the condtion (4.4) is satisfied, this is, however, equivalent to 
$\, 
b\star a = (c\star a)\star a = c\quad ,
\,$
so both relations are equivalent. One can also prove this claim by iterpreting  $\, J_Q(L)\,$ as  $\, Hom_q( \, X(L),\, Q \, )\,$ and using the observations from \cite{K1}, 
Section 3.

Observe that quandles given by symmetric manifolds (Example 2.1) all satisfy the additional condition (4.4).

\bigskip

A topological quandle $\, Q\,$ acts on a product $\, \prod\limits ^n  Q\,$ (any $n$) by $\, (x_1, x_2, ..., x_n)\star a = (x_1\star a, x_2\star a, ..., x_n\star a)\,$. 
(See definition 2.2.) If a link $\, L\,$ can be represented by a closure of a braid on $n$ strands  then the space $\, J_Q(L)\,$ is a subspace of  $\, \prod\limits ^n  Q\,$.

\begin{prop}
The action of $\, Q\,$ on  $\, \prod\limits ^n  Q\,$ restricts to an action of  $\, Q\,$ on  $\, J_Q(L)\,$.
\end{prop} 

\begin{pf}
If $\, \sigma _i\in B_n,\,\, i = 1,..., n-1\,$, is an elementary braid, $\, x\in \prod\limits ^n  Q\,$ and $\, a\in Q\,$ then 
\begin{equation*}
\begin{split}
\sigma _i(\,  x\star a\, )& = \sigma _i(x_1\star a, ... ,x_n\star a )=\\
&=(x_1\star a,..., x_{i-1}  \star a,  x_{i+1}  \star a,  (x_i\star a)\star ( x_{i+1} \star a), x_{i+2} \star a,..., x_n \star a)=\\
&=(x_1\star a,..., x_{i-1}  \star a,  x_{i+1}  \star a,  (x_i\star  x_{i+1}) \star a, x_{i+2} \star a,..., x_n \star a)=\\
&=(x_1,..., x_{i-1},  x_{i+1},  x_i\star  x_{i+1}, x_{i+2},..., x_n ) \star a =\\
& = \sigma _i(\,  x\, ) \star a \qquad .
\end{split}
\end{equation*}
It follows that, for any braid  $\,\, \sigma \in B_n\,$, $\, x\in \prod\limits ^n  Q\,$ and $\, a\in Q\,$ we have 
\begin{displaymath}
 \sigma (x\star a) = \sigma (x)\star a \quad .
\end{displaymath}

If the link $\, L\,$ is the closure of a braid  $\,\, \sigma \in B_n\,$ and if $\, x\in J_Q(L) \subset  \prod\limits ^n  Q\,$ then $\, \sigma (x) = x\,$ and 
\begin{displaymath}
\sigma (x\star a) = \sigma (x)\star a = x \star a
\end{displaymath}
for all $\, a\in Q\,$.
Thus $\,  x \star a\,$ is a fix-point of  $\sigma \,$ and, hence,  $\,  x \star a \in J_Q(L)  \,$ as claimed. 

The action of $\, Q\,$ can also be seen directly if we interpret $\, J_Q(L)  \,$ as the space of colourings of a diagram of $\, L\,$. Given such a colouring, 
an element $\, a\in Q\,$ acts on it by applying the operation $\, ( \,\, )\star a \,$ to all colours.  
\end{pf}

\medskip

Let $\, L_1\,$ and  $\, L_2\,$ be two links in $\, \R ^3\,$ and let $\, L_1\amalg L_2\,$ be their disjoint sum.

\smallskip
\begin{prop}
The invariant space $\, J_Q( L_1\amalg L_2) \,$ of the disjoint sum  $\, L_1\amalg L_2\,$ is homeomorphic to the product $\,  J_Q(L_1) \times  J_Q(L_2)\,$ of 
the invariant spaces of $\, L_1\,$ and $\, L_2\,$. 
\end{prop}
\begin{pf}
If the link $\, L_1\,$ is a closure of a braid $\, \gamma _1\in B_n\,$ on $n$ strands and the link  $\, L_2\,$ is a closure of a braid $\, \gamma _2\in B_m\,$ on $m$ 
strands then  $\, L_1\amalg L_2\,$ is the closure of the braid $\, \gamma _1 \amalg \gamma _2 \in B_{n+m}\,$ om $n+m$ strands. It is immediate from the definition 
of the action of $\, B_{n+m}\,$ on $\, \prod\limits ^{n+m} \, Q\,$ that the space of fix-points of  $\, \gamma _1 \amalg \gamma _2 \,$ is the product of fix-point spaces 
of $\, \gamma _1 \,$ and $\, \gamma _2 \,$.
\end{pf}

\bigskip
 The invariant spaces   $\,  J_Q(L)  \,$ allow one to  construct some numerical (algebraic) invariants of oriented links. We can, for example, 
consider the Poincar\'{e} 
polynomial of  $\,  J_Q(L)  \,$ with coefficients in a field $\, k\,$:
\begin{displaymath}
P_Q( L, \,  k )(t) = \sum\limits_{n=0}^{\infty} \,\, (\dim _kH^n(\, J_Q(L)\, , k )) \, t^n  \qquad .  
\end{displaymath}
It is a polynomial invariant of oriented links.
As follows from Proposition 4.5 we have 
\begin{displaymath}
  P_Q( L_1 \amalg L_2, \,  k )(t) = P_Q( L_1, \,  k )(t) \cdot  P_Q( L_2, \,  k )(t) 
\end{displaymath} 
for a disjoint sum $\, L_1 \amalg L_2 \,$ 
of two links.
\bigskip

\noindent
{\bf Remark 4.4 :}
A special case of the invariant space $\,  J_Q(L)  \,$ of a link  was considered before by D.Silver and S.Williams in \cite{S1}. It was not formulated 
in terms of quandles. When translated into the context of our paper,  the invariant spaces (the spaces of colourings) of Silver and Williams come from 
topological quandles of the type described in Example 2.3 (Alexander quandles) with the topological group $\, G\,$ being the infinite countable product of circles, 
$\, G=\prod\limits _{-\infty}^{\infty} \, T\,$, and the automorphism $\, \sigma : G \longrightarrow G \,$ being the shift $(\sigma (t))_i = t_{i-1}\,$. 
Silver and Williams consider also spaces of colourings of links coming from ``periodic subquandles'' of $\, (G, \sigma )\,$ defined as follows.    
Given a positive integer $\, r\,$, define the $r${\it -periodic subquandle} $\, (G, \sigma )_r\,$ of  $\, (G, \sigma )\,$ to consist of those elements 
$\, t= (t_i)_{-\infty}^{\infty}\,$  from  $\, G =  \prod\limits _{-\infty}^{\infty} \, T\,$ which satisfy 
\begin{displaymath}
t_{i+r} = t_i \qquad \text{and }\,\,\qquad
\sum_{j=i}^{i+r-1}\,  t_j = 0 \qquad \qquad \text{for all }\,\, -\infty < i < \infty \quad .
\end{displaymath} 
The last summation is taken in the circle group $\, T\,$ written additively. 

Silver and Williams define also the spaces of colourings of links with quandles similar to $\, (G, \sigma )\,$ and   $\, (G, \sigma )_r\,$    as above but with 
the circle group $\, T\,$ exchanged for an arbitrary topological group $\, H\,$. All this corresponds to quandles of the type considered in Example 2.3.

\medskip

\begin{lma} 
If $\, G\,$ is an arbitrary topological group and $\, Q_G\,$ is the conjugation quandle of $\, G\,$ then, for any link $\, L\,$ in $\, \R ^3\,$, the invariant space 
$\, J_{Q_G}(L)\,$ is homeomorphic to the space of group homomorphisms $\, Hom(\,  \pi _1(\R ^3 - L)\, , \, G\, )\,$ of the fundamental group of the link complement 
$\,\R ^3 - L\,$  into $\, G\,$. 
\end{lma}
\begin{pf}
According to Corollary 4.1 the invariant space $\, J_{Q_G}(L)\,$ is homeomorphic to the space of quandle homomorphisms $\, Hom_q(\, X(L)\, ,\, Q_G \, ) \,$. If $\, D\,$ 
is an oriented plane diagram of the link $\, L\,$ and $\, E\,$ is the set of arcs of $\,D\,$  then every element of $\, Hom_q(\, X(L)\, ,\, Q_G \, ) \,$ is given by 
  a map of sets 
$\, \varphi :E\longrightarrow G\,$ satisfying the relations (4.2)
\begin{displaymath}
\varphi (g) =\varphi (f)\star \varphi (e) = \varphi (e)^{-1}\, \varphi (f) \, \varphi (e) 
\end{displaymath}
at every crossing 
\begin{displaymath}
\xymatrix{
e& g\\
f\ar@{-}[ru]|\hole  & e\ar[lu]
}
\end{displaymath}
\centerline{Figure 3}

\medskip

\noindent
of $\, D\,$  ($\, e,f,g\,$ being arcs of the diagram).

On the other hand the fundamental group of the link complement $\,  \pi _1(\R ^3 - L)   \,$ is a group generated by the set  $\, E\,$ of arcs of the diagram subject to 
relations $\, e\, f = g\, e\,$ for any crossing of $\, D\,$ as in Figure 3. Hence every element of   $\, Hom(\,  \pi _1(\R ^3 - L)\, , \, G\, )\,$ is given by 
a map of sets $\, \psi :E\longrightarrow G\,$ satisfying the relations $\, \psi (g) = \psi (e)\, \psi (f)\, \psi (e)^{-1}\,$. 
Hence , taking $\, \varphi (h) = \psi (h)^{-1}\,$ 
for any  $\, h\in E\,$ we get $\, \varphi \in Hom_q(\, X(L)\, ,\, Q_G \, ) \,$. This gives a homeomorphism of  $\, Hom(\,  \pi _1(\R ^3 - L)\, , \, G\, )\,$ with       
$\, Hom_q(\, X(L)\, ,\, Q_G \, ) \,$ and, by Corollary 4.1, with $\, J_{Q_G}(L)\,$.
\end{pf}

\bigskip

\section{Examples}

\bigskip   
{\bf Example 5.1 :}  Let $\, O_n\,$ be the trivial link with $n$ components. It is the closure of the trivial (identity) braid on $n$ strands, $\, i_n\in B_n\,$. 
For any quandle $\, Q\,$
the braid $\, i_n\,$  acts as the identity on the product $\, \prod\limits ^n Q \,$. Hence, for any quandle  $\, Q\,$ the invariant space  $\, J_Q(O_n)\,$, being 
the space of fixed points of this map, is equal to  $\, \prod\limits ^n Q \,$.     

\medskip

In this section we take the topological quandle $\, Q\,$ to be the 2-dimensional unit sphere $\, S^2\,$ in $\, \R ^3\,$ with its structure of the symmetric Riemannian 
manifold  inducing the quandle structure (as in Example 2.1). We shall determine the invariant spaces $\, J_Q(L)\,$ for the Hopf link and for two knots: 
the trefoil knot and 
the figure 8 knot. Since the quandle  $\, Q = S^2\,$ is a {\it kei} i.e. it satisfies the additional condition (4.4), the invariant spaces  $\, J_{S^2}(L)\,$ are independent 
of the orientation of the link $\, L\,$ (see Remark 4.3).

\bigskip

{\bf Example 5.2 :} Let $\, H\,$ be the Hopf link (see Figure 4). 

\bigskip

\medskip
$$\includegraphics[width=3.7cm, height=1.8cm]{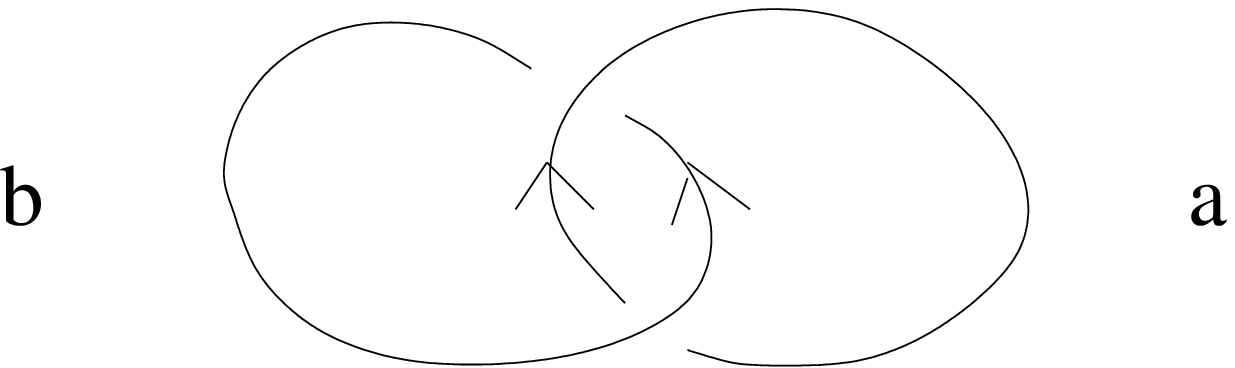}$$

\centerline{Figure 4}

\bigskip

\begin{prop} 
The invariant space  $\, J_{S^2}(H)\,$ of the Hopf link is homeomorphic to the disjoint union $\, S^2 \amalg S^2\,$ of two copies of $\, S^2\,$.
\end{prop}

\begin{pf}
We interpret the space $\, J_{S^2}(H)\,$ as the space of colourings of the diagram of the link with colours from $\, S^2\,$.  The diagram of the Hopf link has two arcs. 
We assign to them elements $\, a, b\in S^2\,$, respectively, as in Figure 4. There are two crossings. The upper crossing 
 gives an 
equation $\, b\star a = b\,$. The lower one gives  an equation $\, a\star b = a\,$. The only solutions to each one of this equations are $\, b = \pm a\,$ for arbitrary 
$\, a\in S^2\,$. Hence the space   $\, J_{S^2}(H) = \{ (a,\, \pm a)\in S^2 \times  S^2 \, | \, a \in S^2 \, \} = S^2 \amalg S^2\,$.
\end{pf}

\medskip
{\bf Example 5.3 :} The trefoil knot $\, K_3\,$ is a closure of the braid $\, \sigma = \sigma _1 ^{-3} \in B_2\,$ on two strands:

\bigskip
$$\includegraphics[width=2.5cm, height=3.5cm]{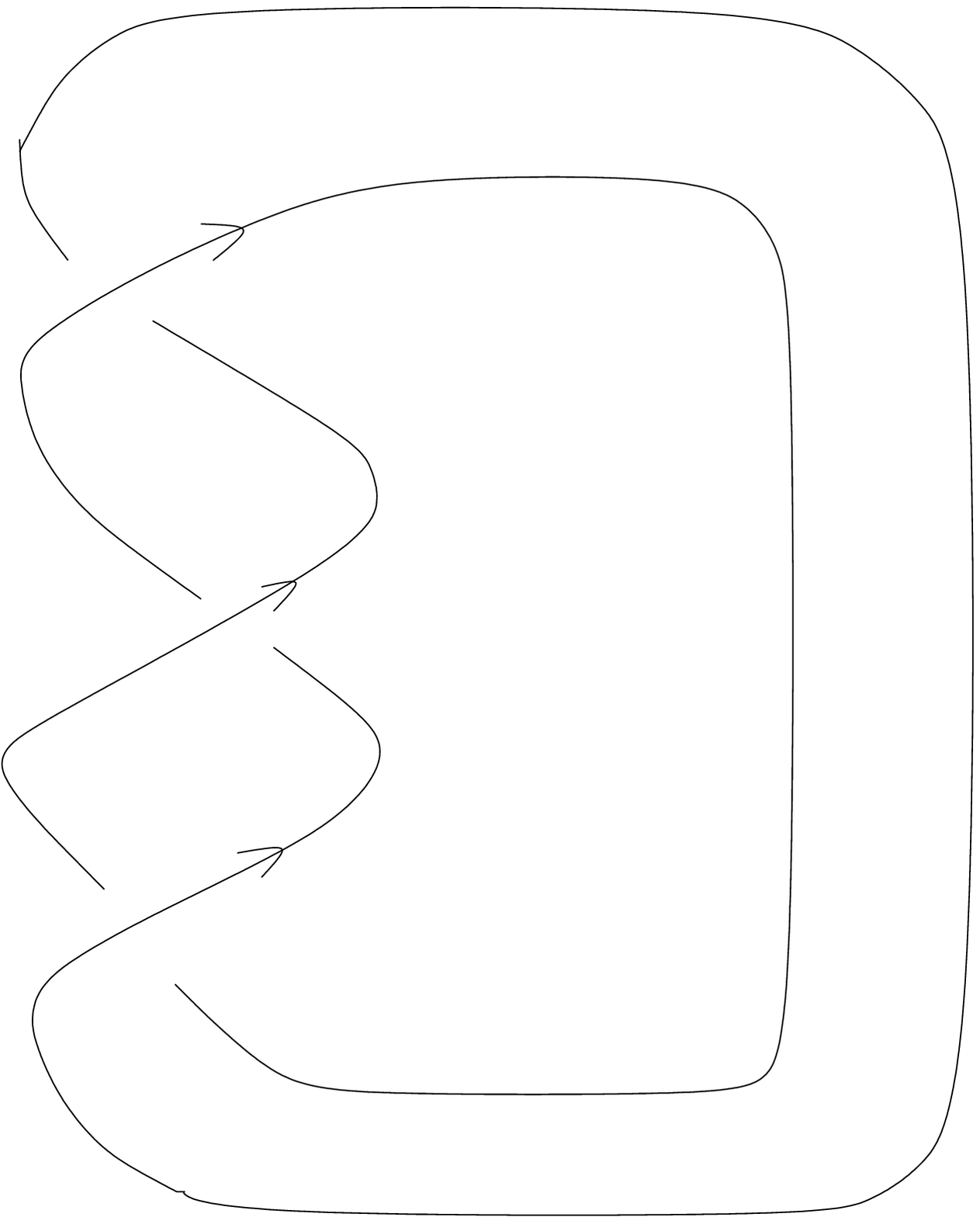}$$

\medskip

\centerline{Figure 5}

\bigskip

\begin{prop}
The invariant space  $\, J_{S^2}(K_3)\,$ of the trefoil knot  has two connected components, one diffeomorphic to $\, S^2\,$ and the other to   $\, \R P^3\,$.
\end{prop}

\begin{pf}
The space $\, J_{S^2}(K_3)\,$ is the subspace of fix-points of the braid $\,\sigma = \sigma _1 ^{-3}\,$ acting on $\, Q\times Q\,$. That is the same as 
the subspace of fix-points  of the braid   $\, \sigma _1 ^3\,$. For $\, (a, b) \in Q\times Q\,$ we have
\begin{equation*}
\begin{split}
\sigma _1^3 (a, b) &= \sigma _1^2 (b, a\star b) = \sigma _1 (a\star b,\, b\star ( a\star b)) =\\ &= (\,  b\star ( a\star b), \, (a\star b)\star ( b\star ( a\star b)) \, )
\end{split}
\end{equation*}
Therefore  $\,  (a, b) \in Q\times Q\,$ is a fix-point of $\, \sigma\,$ if and only if 
\begin{equation}
\left\{\begin{array}{l}
 b\star ( a\star b)= a\\
 (a\star b)\star a = b
\end{array}\right.
\end{equation}
If $\, b = a\,$ then $\, a\star b = a\,$ and both equations of (5.1) are satisfied. If $\, b=-a\,$ then $\, a\star b =a\,$ and $\,  (a\star b)\star a = a\star a = a \ne b\,$ 
and the equations (5.1) are not satisfied. Finally, if $\, b\ne \pm a\,$ then there is a unique great circle of the sphere $\, S^2\,$ containing $\, a\,$ and $\, b\,$. 
That circle contains also $\, a\star b, \, (a\star b)\star a\,$ and $\,  b\star ( a\star b)\,$. Let us parametrise the circle with the angle starting at the point $\, b\,$. 

\bigskip
$$\includegraphics[width=7cm, height=4cm]{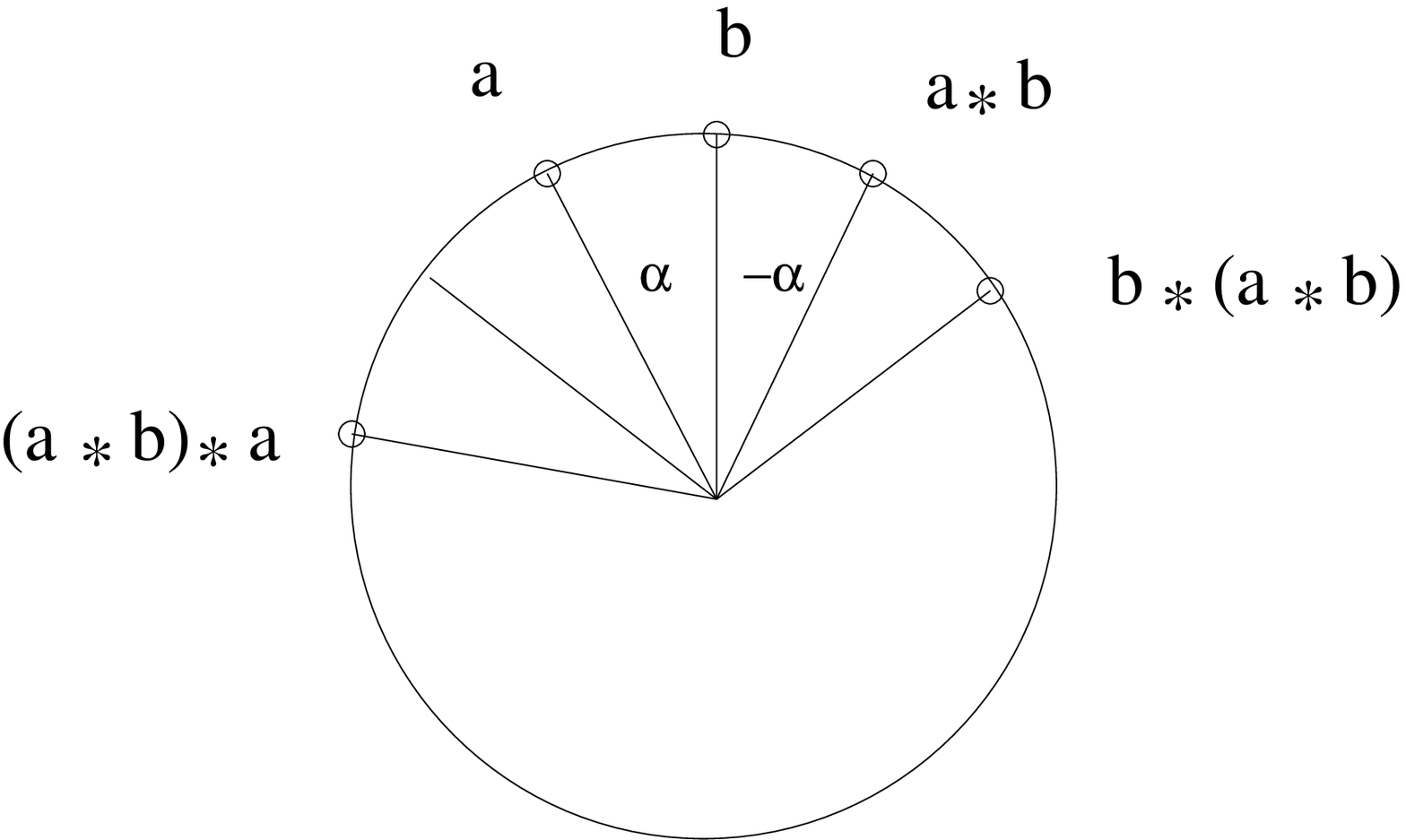}$$
\centerline{Figure 6 }

\bigskip

If  $\, a\,$ has an angle $\,\alpha\,$ then $\, a\star b\,$ has the  angle  $\, -\alpha\,$ , the point $\,  b\star ( a\star b)\,$ has the angle $\, -2\alpha\,$ and 
$\,  (a\star b)\star a\,$ has the angle  $\, 3\alpha\,$. (See Figure 6.) 
Hence the equations (5.1) are satisfied if and only if $\, -2\alpha \equiv \alpha\,$ and $\, 3\alpha \equiv 0\,$ 
(mod $2\pi$) that is  iff $\, \alpha = 0 \,\,\,\text{or}\,\, \alpha = \pm 2\pi /3\,$. 

Let $\, p: J_Q(K_3) \longrightarrow Q\,$ be the projection on, say, the first coordinate in $\, Q\times Q\,$. Given $\, a\in Q\,$, it follows from the calculations above 
that the fibre $\, p^{-1}(a) \subset Q\,$ is the union the point $\, a\,$ and the circle laying on $\, Q=S^2\,$ in a plane perpendicular to $\, a\,$. The plane has 
distance $\, \frac 12\,$ to the origin and  $\, \frac 32\,$ to $\, a\,$. The union of all the circles (over all $\, a\in Q\,$) is diffeomorphic to the unit circle 
subbundle of the tangent bundle of $\, S^2\,$ and hence diffeomorphic to $\, \R P^3\,$. That is the component of $\, J_Q(K_3)\,$ consisting of solutions 
$\, (a,b)\in Q\times Q\,$ satisfying $\, a\ne b\,$. The other component of $\, J_Q(K_3)\,$ consists of solutions of the form $\, (a,a)\in Q\times Q\,$ 
and is diffeomorphic to $\, Q=S^2\,$. 
\end{pf}

\bigskip

{\bf Example 5.3 :} Let $\, K_{(-2.2)}\,$ be the figure-eight knot  (see Figure 7). It is the two-bridge knot with the Conway normal form $\, C(-2,2)\,$ 
(see \cite{K3}). 
It has a diagram $\, D\,$ with 4 arcs and 4 crossings

\bigskip
$$\includegraphics[width=3.5cm, height=5cm]{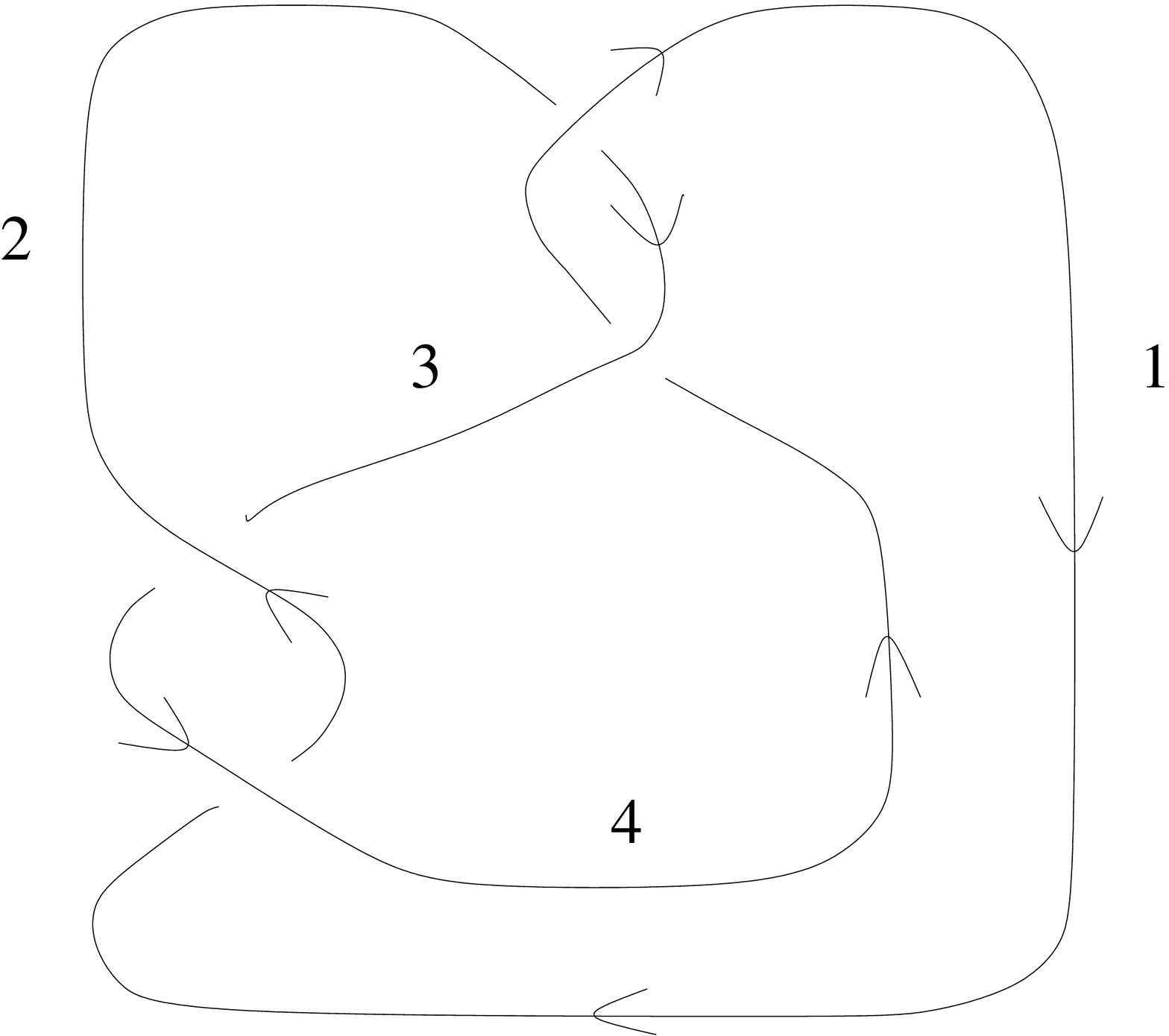}$$

\medskip

\centerline{Figure 7 }

\medskip

\begin{prop}
The invariant space  $\, J_{S^2}(K_{(-2.2)})\,$ of the figure-eight  knot   $\, K_{(-2.2)}\,$  has three connected components, 
one diffeomorphic to $\, S^2\,$ and each one of 
the other 
two diffeomorphic to   $\, \R P^3\,$.
\end{prop}

\begin{pf} We interpret  $\, J_{S^2}(K_{(-2.2)})\,$ as the space of colourings of the diagram  of  $\, K_{(-2.2)}\,$.
We number the arcs of the diagram as indicated in Figure 7 and we number the crossings from the highest to the lowest.  If we assign a colour 
$\, a\in Q\,$ to the arc  1 and  
a colour $\, b\in Q\,$ to the arc 2 then the first (highest) crossing forces us to assign the colour $\, b\star a\,$ to the arc 3. If we assign a colour $\, c\in Q\,$ 
to the arc 4  then the second  crossing gives us an equation $\, c\star (b\star a)= a\,$. That determines $\, c\,$ uniquely.  
The third crossing gives us an equation $\, c\star b = b\star a\,$ and the 
fourth crossing the equation $\, b\star c = a\,$. Hence we have three colours $\, a, b, c\,$ which have to satisfy the system of equations
\begin{equation}
\left\{\begin{array}{l}
c\star (b\star a) = a\\
b\star c =a\\
c\star b = b\star a
\end{array}\right.
\end{equation}

However, in this quandle $\, Q\,$ 
the operation $\, (\,\, )\star d\,$ is involutive for every $\, d\in Q\,$ i.e. $\, (e\star d)\star d = e\,$ for all $\, e,d\in Q\,$. Therefore, applying the operation 
$\, ( \,\, \, )\star (b\star a)\,$ to both sides of the first equation we get that the colour assigned to the arc 4 is $\, c=a\star (b\star a)\,$. 
We are then left with two colours $\, a\,$ and $\, b\,$ and two equations
\begin{equation}
\left\{\begin{array}{l}
b\star (a\star (b\star a)) =a\\
(a\star (b\star a))\star b = b\star a.
\end{array}\right.
\end{equation} 

We solve equations (5.3) in a way similar to how we solved equations (5.1).        

If $\, b = a \,$ then equations (5.2) are satisfied due to $\, a\star a = a\,$. If $\, b= -a \,$ then $\, b\star a=b, \,\, a\star (b\star a)=a\star b=a\,$ and 
$\, b\star ( a\star (b\star a))=b\star a = b = -a \ne a\,$. Thus the first equation of (5.1) is not satisfied if $\, b=-a\,$.

If $\, b\ne \pm a\,$ then $\, a\,$ and $\, b\,$ belong to a unique great circle of $\, Q=S^2\,$. The circle contains also all 
the expressions in both sides of equations (5.3). We parametrise the circle with the angle, the point $\, a\,$ 
corresponding to angle $0$. 

\bigskip

$$\includegraphics[width=7cm, height=4cm]{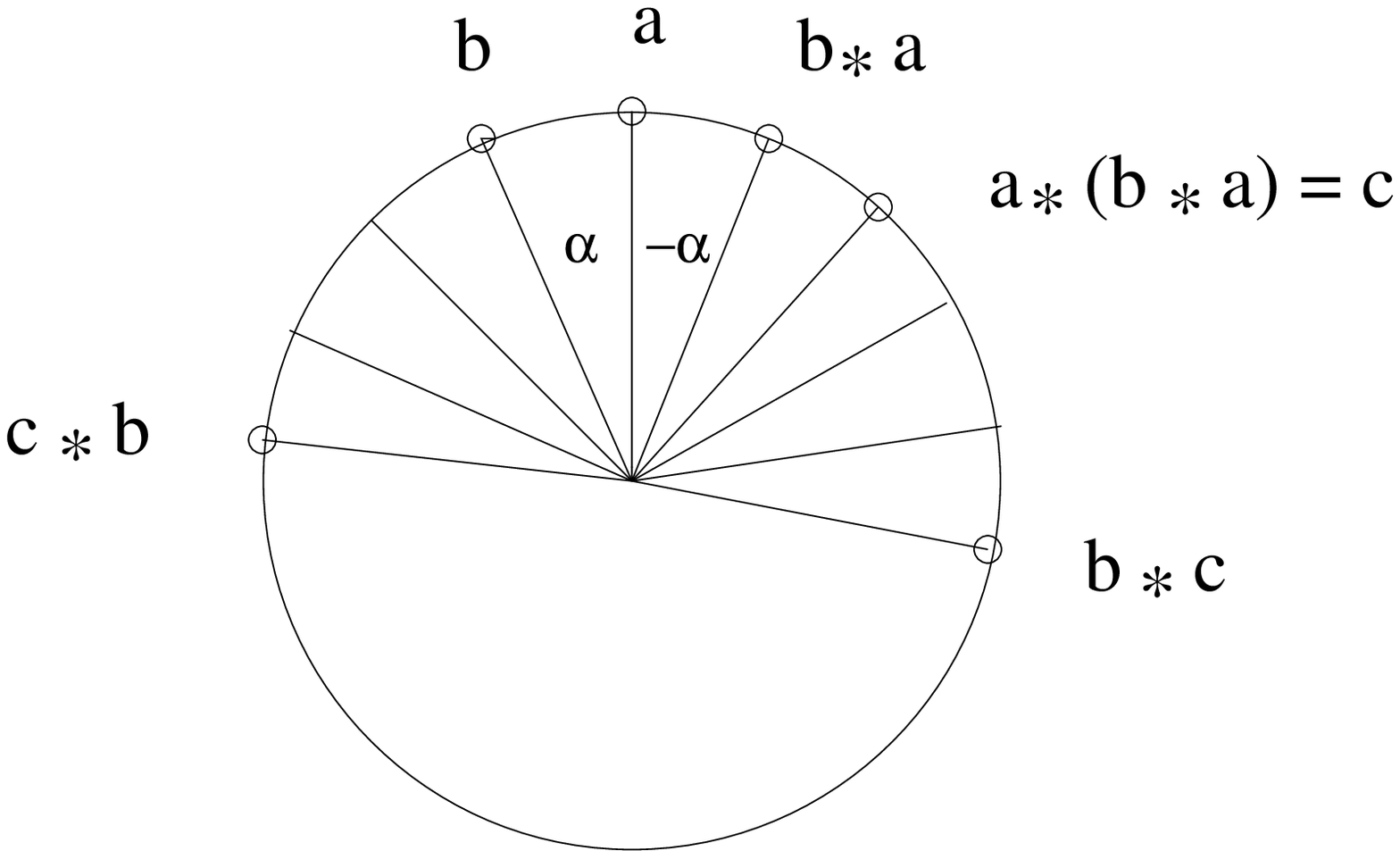}$$

\medskip

\centerline{Figure 8}

\medskip

If $\, b\,$ has angle $\, \alpha \,$ then $\, b\star a \,$ has angle $\, -\alpha \,$, $\,  a\star (b\star a)\, $ has angle $\, -2\alpha \,$,  
$\, b\star ( a\star (b\star a))\,$ has angle $\, -5\alpha \,$ and the point $\, ( a\star (b\star a))\star b \,$ has angle $\, 4\alpha \,$. (See Figure 8.)
Hence equations (5.3) are satisfied 
if and only if $\, -5\alpha \equiv 0\,$ and $\, 4\alpha \equiv -\alpha \,$ (mod $2\pi$) i.e. iff $\, \alpha = 2\pi k/5\, , \,\, k$ being an integer.

Thus, given any $\, a\in Q=S^2\,$, there are 5 points $\, b\in Q\,$ on every great cicle of $\, S^2\,$ through $\, a\,$ such that the pairs $\, (a,b)\,$ solve 
equations (5.3). These 5 points $\, b\,$ are vertices of a regular pentagon, one of them being equal to $\, a\,$. Hence, given any $\, a\in Q=S^2\,$ 
the set of points $\, b\in Q\,$ such that $\, (a,b)\,$ solve 
equations (5.3) is a disjoint union of two circles and of  the point $\, a\,$.

The invariant space $\, J_Q(K_{(-2,2)})\,$ is a subspace of $\, Q\times Q\,$. If $\, p:  J_Q(K_{(-2,2)}) \longrightarrow Q\,$ is the projection on the first coordinate in 
 $\, Q\times Q\,$ then, for any   $\, a\in Q\,$,  the fibre $\, p^{-1}(a)\,$ is the disjoint union of $\, a\,$ and of two circles, as described above. 

Therefore, the space 
 $\, J_Q(K_{(-2,2)})\,$ has three connected components: one diffeomorphic to $\, Q=S^2\,$ (equal to the diagonal of  $\, Q\times Q\,$) and two other components, each one 
being diffeomorphic to the unit circle subbundle  of the tangent bundle of $\, S^2\,$, i.e. each one being diffeomorphic to $\, \R P^3\,$.     
\end{pf} 

\medskip

 \noindent
{\bf Remark 5.1 :} The trefoil knot and the figure-eight knot are examples of two-bridge knots. There is a rather simple algorithm for computation of the invariant space 
$\, J_{S^2}(L)\,$ for all two bridge knots and links from their diagrams in the Conway normal form. The algorithm follows the pattern we have used in computations 
for the trefoil and the figure-eight knot. The result is that the space  $\, J_{S^2}(L)\,$ is a disjoint union of one or two copies of $\, S^2\,$ and a number 
of copies of   $\, \R P^3\,$.  

\bigskip

\section{A further example}

\bigskip

In this section we consider a topological quandle $\, Q\,$ which is a conjugacy class in the Lie group $\, SL(2,\, \C )\,$ and we compute the invariant space 
 $\, J_Q(K_3)\,$ for the trefoil knot $\, K_3\,$. 

\medskip

Let $\, C\,$ be the matrix $\, \bigl( \begin{smallmatrix} 1&0\\ 1&1 \end{smallmatrix} \bigr) \in  SL(2,\, \C )\,$ and let $\, Q\,$ be the conjugacy class of $\,C\,$ in  
$\, SL(2,\, \C )\,$. The space  $\, Q\,$ is a topological quandle with the quandle operation $\, h\star g: = g^{-1}hg\,$ for $\, h,g\in Q\,$ (see Example 2.5).

\medskip
\noindent
For each $\, \bigl( \begin{smallmatrix} a&b\\ c&d \end{smallmatrix} \bigr) \in  SL(2,\, \C )\,$ we have
\begin{equation}
\begin{split}
\left( \begin{array}{cc} a&b \\ c&d\end{array}\right)^{-1}\left( \begin{array}{cc} 1&0 \\ 1&1\end{array}\right)\left( \begin{array}{cc} a&b \\ c&d\end{array}\right) &= 
\left( \begin{array}{cc} d&-b \\ -c&a\end{array}\right)\left( \begin{array}{cc} 1&0 \\ 1&1\end{array}\right)\left( \begin{array}{cc} a&b \\ c&d\end{array}\right) =\\
%&=\left( \begin{array}{cc} d-b&-b \\ a-c&a\end{array}\right)\left( \begin{array}{cc} a&b \\ c&d\end{array}\right)=\\
%&=\left( \begin{array}{cc} ad-ab-bc&db-b^2-bd \\ a^2-ac+ac&ab-bc+ad\end{array}\right)=\\
&=\left( \begin{array}{cc} 1-ab&-b^2 \\ a^2&1+ab\end{array}\right)
\end{split}
\end{equation}
and the entries of the last matrix determine the pair of complex numbers $\, (a,b)\,$ uniquely up to simultaneous multiplication by $\, \pm 1\,$.

It follows that every point in $\, Q\,$ can be identified with an ordered pair of complex numbers $\, (a,b)\,$ well-defined up to multiplication by $\, \pm 1\,$, 
subject only 
to the condition   $\, (a,b)\ne (0,0)\,$.

Let us  consider the space $\, \C ^2-\{\vec 0 \}\,$ with the equivalence relation $\, (a,b)\sim (e,f) \,\, \Leftrightarrow \,\, (a,b)=\pm (e,f)\,$.
We have the map $\,\varphi : Q \longrightarrow  (\C ^2-\{\vec 0 \})/\sim \,$, 
\begin{displaymath} 
\varphi \bigl(\bigl( \begin{smallmatrix} 1-ab&-b^2\\ a^2&1+ab \end{smallmatrix} \bigr)\bigr) = [\, a,b\, ] \quad.
\end{displaymath}

\begin{prop} The map $\,\varphi : Q \longrightarrow  (\C ^2-\{\vec 0\})/\sim \,$ is a homeomorphism.
% The quandle $\, Q\,$ is homeomorphic to the quotient space $\, \C ^2-\{\vec 0\}/\sim \,$. 
The quandle operation on $\, Q\,$ corresponds to 
\begin{displaymath}
 [\, a,b\, ]\star [\, c,d\, ] = [\, a-acd+bc^2\, , \, b-ad^2+bcd\, ]
\end{displaymath}
on $\, (\C ^2-\{\vec 0 \})/\sim \,$.
\end{prop}

\begin{pf}
The first claim has already been proven. The second one follows from 
\begin{displaymath}
\begin{split}
 \varphi ^{-1}([\, a,b\,])\, \star&\,  \varphi ^{-1}([\, c,d\,])=   \left( \begin{array}{cc} 1-ab&-b^2 \\ a^2&1+ab\end{array}\right)\, \star \, 
\left( \begin{array}{cc} 1-cd&-d^2 \\ c^2&1+cd\end{array}\right) =     \\
 &=\left( \begin{array}{cc} 1-cd&-d^2 \\ c^2&1+cd\end{array}\right) ^{-1} \left( \begin{array}{cc} 1-ab&-b^2 \\ a^2&1+ab\end{array}\right) 
\left( \begin{array}{cc} 1-cd&-d^2 \\ c^2&1+cd\end{array}\right) =\\ &= \left( \begin{array}{cc} 1-ef&-f^2 \\ e^2&1+ef\end{array}\right)=\\
&= \varphi ^{-1}([\, e,f\,])
\end{split}
\end{displaymath} 
with $\, e=a-acd+bc^2\,$ and $\, f=b-ad^2+bcd\,$. Only the fourth identity has to be checked, others are the matter of definitions. That is left to the reader. 

%The second one is obtained by conjugating the matrix  
%$\, \bigl( \begin{smallmatrix} 1-ab&-b^2\\ a^2&1+ab \end{smallmatrix} \bigr)\,$ 
%by the matrix $\, \bigl( \begin{smallmatrix} 1-cd&-d^2\\ c^2&1+cd \end{smallmatrix} \bigr)\,$. The details are left to the reader. 
\end{pf}

\medskip

\noindent
{\bf Remark 6.1 :} 
%\noindent
%{\bf 1.}
The manifold  $\, (\C ^2-\{\vec 0 \})/\sim \,$ is diffeomorphic to the product $\, \R P^3 \times \R\,$.

%\noindent
%{\bf 2.}  $\, (\C ^2-\{\vec 0 \})/\sim \,$ has a structure of a complex manifold. The quandle operation is a holomorphic map.

\bigskip

We shall now describe the invariant space $\, J_Q(K_3)\,$ of the trefoil knot $\, K_3\,$.

\medskip
 
\begin{prop}
The invariant space  $\, J_Q(K_3)\,$ of the trefoil knot  $\, K_3\,$ has two connected components , one diffeomorphic to  $\, \R P^3 \times \R\,$ and 
the other one diffeomorphic to  
$\, \R P^3 \times \R ^3\,$. 
\end{prop}

\begin{pf}
We refer again to Figure 5. The begining of the proof is the same as that of Proposition 5.1. The trefoil knot  $\, K_3\,$ is the closure of the braid 
$\, \sigma = \sigma _1^{-3}\in B_2\,$ on two strands. Hence the space  $\, J_Q(K_3)\,$  is the subspace of fix-points of the braid $\, \sigma _1^{-3}\,$ 
acting on $\, Q\times Q\,$. That is the same as the subspace of fix-points of the braid  $\, \sigma _1^3\,$. As in the proof of Proposition 5.1 we see that it consists 
of points $\, (g,h)\in  Q\times Q\,$ satisfying 
\begin{equation}
\left\{\begin{array}{l}
 h\star ( g\star h)= g\\
 (g\star h)\star g = h.
\end{array}\right.
\end{equation}
(Compare (5.1).) Since, in the quandle $\, Q\subset SL(2,\C)\,$, one has $\, g\star h = h^{-1}gh\,$ with the product being the matrix multiplication, we see that 
each of the equations of (6.2) is equivalent to $\, hgh=ghg\,$. Thus both equations of (6.2) are equivalent and the system (6.2) is equivalent to
\begin{equation}
(g\star h)\star g = h \qquad .
\end{equation} 
For any $\, g\in Q\,$ the pair $\, (g,h) = (g,g)\in Q\times Q \,$ is a solution to (6.3). We shall now look for the remaining solutions.

Let us first observe that if $\, \psi :Q\longrightarrow Q\,$ is an automorphism of the quandle $\, Q\,$ and $\, (g,h)\in  Q\times Q \,$ is a solution to (6.3) then so is 
even the pair $\, (\psi (g), \, \psi (h))\,$. (Compare with Remark 4.1.) Secondly, the group $\, SL(2,\C)\,$ acts by conjugation on the quandle $\, Q\,$ as a group of
its quandle  automorphisms. The action is transitive by the definition of $\, Q\,$. Therefore, if we find all solutions to (6.3) with one fixed $\,g = g_0\,$, we shall get 
all solutions (for all $\, g\,$) through the action of  $\, SL(2,\C)\,$.

We now denote by  $\, Y\,$ the space   $\, (\C ^2-\{\vec 0 \})/\sim \,\,\,$ and we identify the quandle $\, Q\,$ with the space  $\, Y\,$ equipped with 
the quandle operation as described in Proposition 6.1.

Let us choose $\, g_0= [1,0]\in Y \,$. Let $\, h=[\alpha \,  ,\, \beta ]\in Y\,$. Then, 
according to Proposition 6.1,  $\, g_0\star h = [1,0]\star [\alpha \,  ,\, \beta ]  = [ 1-\alpha \beta \, ,\,  -\beta ^2 \, ]\,$ and 
$\, (g_0 \star h) \star g_0  = [ 1-\alpha \beta \, ,\,  -\beta ^2 \, ] \star [1,0] =   [1- \alpha \beta -\beta ^2 \, , \, - \beta ^2\, ]\,$. 
Thus the pair $\, (g_0,h)\,$ satisfies the equation 
(6.3) 
if and only if 
\begin{equation}
 [1- \alpha \beta - \beta ^2 \,  , \,  - \beta ^2 ] = [\alpha \, , \, \beta ] \qquad .
\end{equation} 
The second coordinate implies $\, \pm \beta ^2 =\beta \,$. Hence  $\, \beta = 0,\pm 1\,$. If $\, \beta =0\,$, the first coordinate gives $\, \alpha =\pm 1\,$. We get 
$\, h=[1,0]=g_0\,$ and $\, (g_0,h)=(g_0, g_0)\,$, which is the solution we have discussed earlier. If $\, \beta = \pm 1\,$, then (6.4) is satisfied by 
any $\, \alpha \in \C\,$. Thus all solutions $\, (g,h)\in Q\times Q\,$ to (6.3) with $\, g=g_0=[1,0]\,$ are given by $\, h=g_0=[1,0]\,$ or by $\, h=[\alpha \, ,\, 1]\,$ with 
arbitrary $\, \alpha \in \C\,$. 

For any pair of complex numbers  $\, (a,b)\in  \C ^2-\{\vec 0 \}\,$,\, let $\, {\tilde a} = \dfrac{\overline a}{|a|^2+|b|^2}\,$, 
let  $\, {\tilde b} = \dfrac{\overline b}{|a|^2+|b|^2}\,$ and let $\, A_{(a,b)}\,$ denote the matrix
\begin{displaymath}
 A_{(a,b)} = \left( \begin{array}{cc}
a&b \\
-{\tilde b} &{\tilde a}
\end{array}\right)
\in SL(2,\C \, )\qquad .
\end{displaymath}    
  
Under the homeomorphism $\, \varphi : Q \longrightarrow Y\,$ of Proposition 6.1 the point $\, g_0= [1,0]\,$  in $\, Y\,$ corresponds  
to the matrix 
$\, C = \bigl( \begin{smallmatrix} 1&0\\ 1&1 \end{smallmatrix} \bigr)\,$ in $\, Q\,$, while the point $\, h=[\alpha\, , 1 ]\,$ in $\, Y\,$ corresponds to the matrix 
$\, D_{\alpha}=   \bigl( \begin{smallmatrix} \alpha &1\\ -1&0 \end{smallmatrix} \bigr)^{-1}  \bigl( \begin{smallmatrix} 1&0\\ 1&1 \end{smallmatrix} \bigr)
 \bigl( \begin{smallmatrix} \alpha &1\\ -1&0 \end{smallmatrix} \bigr) =  \bigl( \begin{smallmatrix} 1- \alpha &-1\\      \alpha ^2 & 1+\alpha \end{smallmatrix} \bigr) 
\,$ in $\, Q\,$.

The conjugation of $\, C\,$ and of $\, D_{\alpha}\,$     by  $\, A_{(a,b)}\,$ gives, according to (6.1),
\begin{displaymath}
A_{(a,b)} ^{-1}\, C \, A_{(a,b)} =  \left( \begin{array}{cc} 1-ab&-b^2 \\ a^2&1+ab\end{array}\right) = \varphi ^{-1} ([a,b]) 
\end{displaymath}  
and 
\begin{displaymath}
\begin{split}
A_{(a,b)} ^{-1}\,  D_{\alpha}  \, A_{(a,b)}& =  A_{(a,b)} ^{-1}\,  \left( \begin{array}{cc} \alpha & 1 \\ -1&0\end{array}\right)^{-1}\,  
\left( \begin{array}{cc} 1 & 0 \\ 1&1\end{array}\right)\, 
\left( \begin{array}{cc} \alpha  & 1 \\ -1&0\end{array}\right) \, A_{(a,b)} = \\
&=  \left( \begin{array}{cc} \alpha\, a - {\tilde b} & \alpha\, b + {\tilde a} \\ -a&-b\end{array}\right)^{-1}\, 
\left( \begin{array}{cc} 1 & 0 \\ 1&1\end{array}\right)\, 
 \left( \begin{array}{cc} \alpha\, a - {\tilde b} & \alpha\, b + {\tilde a} \\ -a&-b\end{array}\right)      = \\
&= \varphi ^{-1} ([\,  \alpha a - {\tilde b}\,    , \,  \alpha b + {\tilde a}\,  ]) \quad .
\end{split}
\end{displaymath} 
Consequently, all solutions $\, (g,h)\in Y\times Y\,$ to the equation (6.3) are either of the form

(i) $\, (g,g)\,$ with arbitrary $\, g\in Y\,$, 

\noindent
or of the form

(ii) $\, g=[\, a,b\, ]\,$ and $\, h = [\,  \alpha a - {\tilde b}\,    , \,  \alpha b + {\tilde a}\,  ]\,$ with arbitrary $\,  g=[\, a,b\, ]\in Y\,$ 
and 

\qquad  $\, \alpha \in \C \,$.   

\noindent
These two sets of solutions are disjoint.

The solutions of type (i) form one connected component of $\, J_Q(K_3)\,$ diffeomorphic to $\, Q\,$ and, hence, to $\, \R P^3\times \R\,$.

To describe the space of solutions of type (ii) let us define a map
\begin{displaymath}
F: Y\times \C \longrightarrow Y\times Y 
\end{displaymath}  
by $\, F([a,b]\, , \, \alpha ) = ([a,b]\, , \,[\,  \alpha a - {\tilde b}\,    , \,  \alpha b + {\tilde a}\,  ]\,$. The map $\, F\, $ is a diffeomorphism from the manifold 
$\,  Y\times \C \,$ onto the space of solutions of type (ii) to the equation (6.3) in $\, Y\times Y\,$. 

Therefore the invariant space $\, J_Q(K_3)\,$, being the space of solutions to (6.3), has two connected components, one diffeomorphic to   $\, \R P^3\times \R\,$ and 
one diffeomorphic to $\,  Y\times \C \,$, hence, to  $\, \R P^3\times \R ^3\,$.
\end{pf}

\medskip

In Section 4 we have already mentioned that one can obtain algebraic invariants of oriented links, for example, by taking the Poincar\'{e} polynomials of 
the invariant spaces   $\, J_Q(L)\,$. The  example of  $\, J_Q(K_3)\,$, described in Proposition 6.2, suggests 
that another interesting algebraic invariant of links  
can be obtained by taking the   Poincar\'{e} polynomials of the {\bf one-point compactification} of the invariant spaces   $\, J_Q(L)\,$. 

For example, with $\, Q\,$ being the quandle studied in this section,  the  Poincar\'{e} polynomial of the one-point compactification of the space  $\, J_Q(K_3)\,$ 
of the trefoil knot $\,K_3\,$ 
(with coefficients in the field $\, \Q\,$ of 
rational numbers) is 
\begin{displaymath}
P_{\, Q, cpt}\, (K_3\, ,\, \Q\, )(t) =  1+t+t^3+t^4+t^6 \quad .
\end{displaymath}
That follows directly from Proposition 6.2.

\bigskip

\noindent 
{\bf Remark 6.2 :} In many cases topological quandles have a richer structure than just that of a topological space. For example, 
the quandle $\, Q \subset SL(2\, , \, \C\, )\,$, considered in this section, is a complex manifold with the quandle operation being a holomorphic map. 
The conjugation quandles of algebraic groups have a structure of an algebraic variety. The invariant  spaces  $\, J_Q(L)\,$ of links with quandles 
$\, Q\,$  of that kind will be complex manifolds or varieties, perhaps with singularities. One may expect further invariants of links to be derived 
from those additional structures of the spaces   $\, J_Q(L)\,$.

\end{document}